\documentclass[12pt]{article}

\usepackage{latexsym}
\usepackage{amssymb}

\newcommand{\I}{{\mathcal I}_f}
\newcommand{\R}{\mathbb R} 
\newcommand{\C}{\mathbb C}
\newcommand{\D}{\mathbb D}
\newcommand{\T}{\mathbb T}
\newcommand{\Sm}{{\mathcal N}_+}
\newcommand{\z}{\zeta}
\renewcommand{\o}{\omega}

\title{A question by Alexei Aleksandrov \\ and logarithmic determinants}
\author{Mikhail Sodin\thanks
{Supported by the Israel Science Foundation of the Israel Academy
of Sciences and Humanities under Grant No. 37/00-1.}
}

\begin{document}
\date{}

\maketitle

\bigskip\par\noindent{\bf \S1. Analytic functions represented by Schwarz'
integrals.}

\medskip\par\noindent
This note is motivated by the following question:
{\em when is an analytic function $f$ in the unit disc $\D$ 
represented by the Schwarz integral 
$$
f(z) = \int_{\T} \frac{\z+z}{\z-z}\, d\mu (\z)
\eqno (1.1)
$$
of a real measure $\mu$ on the unit circle $\T$?} 
Two classical necessary conditions are:

\smallskip\par\noindent{\bf (i)} (Smirnov). $f\in \Sm$, that is,
$$
\log|f(z)|\le \int_{\T} \log|f(\z)| \,
\hbox{Re} \left(\frac{\z+z}{\z-z}\right) \, dm(\z)\,,
\qquad z\in\D\,,
\eqno (1.2)
$$
where $m$ is the normalized Lebesgue measure on $\T$, $m(\T)=1$.
(We use the same notation $f$ for an analytic 
function in the unit disc and for its non-tangential boundary values 
on the unit circle.)
In fact, Smirnov proved a stronger result: 
$f\in \bigcap_{p<1} H^p$.

\smallskip\par\noindent{\bf (ii)}
(Kolmogorov). $f\in L^{1,\infty}(\T)$, that is, $m_f(t)=O(t^{-1})$,
$t\to\infty$, where $m_f(t) = m\{\z\in\T:\, |f(\z)|\ge t \}$.
More precisely, 
if $f$ is represented in the form (1.1), then
$$
m_f(t) \le \frac{C||\mu||}{t}\,,
\qquad 0<t<\infty\,,
$$
where $C$ is a numerical constant.

\smallskip
Another necessary condition is 
\smallskip\par\noindent{\bf (iii)} (Hruschev and Vinogradov \cite{HV})
if $f$ us represented in the form (1.1), then  
there exists a limit
$$
\lim_{t\to \infty} tm_f(t) = \frac{||\mu_{\hbox{\small sing}}||}{\pi}\,.
$$
A weaker result that $\lim_{t\to \infty} tm_f(t) = 0$ if the
measure $\mu$ in the representation (1.1) is absolutely continuous, is due
to Titchmarsh \cite{T}.

\smallskip
Besides the above conditions, there are the trivial restrictions:
\smallskip\par\noindent{\bf (iv)}
$$
\hbox{Re} f \left( = \frac{d\mu_{\hbox{\small a.c}}}{dm}\right)
\in L^1(\T)\,, 
\qquad \hbox{and} \qquad  
f(0) \left( =\int_\T d\mu \right) \in \mathbb R\,.
$$

\smallskip
In \cite{A1} and \cite{A2} Aleksandrov proved 

\medskip\par\noindent{\bf Theorem A1. }{\em 
The set of conditions 
(i), (ii), and (iv) is sufficient for representation (1.1), with}
$$
||\mu|| \le ||\hbox{Re} f||_1 + C ||f||_{1,\infty}\,,
$$ 
{\em where $||\,.\,||_1=||\,.\,||_{L^1(\T)}$, 
$||f||_{1,\infty} = \sup_{t>0} tm_f(t)$, and $C$ is a positive
numerical constant.}

\smallskip
That is, an analytic function $f$ is represented by the Schwarz integral
(1.1) iff conditions (i), (ii) and (iv) hold.
Furthermore, Aleksandrov suggested that in the sufficiency part condition
(ii) can be weakened and asked 
{\em whether conditions (i), 
$$
\liminf_{t\to\infty} tm_f(t) < \infty\,,
$$ 
and (iv) already guarantee that $f$ is represented by the Schwarz integral 
(1.1).} Another form of this question is (cf. \cite{A3}):
{\em whether there exists a non-constant analytic function $f$ of
Smirnov's class $\Sm$ such that $\hbox{Re} f = 0$ a.e. on $\T$ and} 
$$
\liminf_{t\to\infty} tm_f(t) = 0\,.
\eqno (1.4)
$$ 
If such a function exists, then in view of the necessary conditions cited
above, it cannot be represented by the Schwarz integral (1.1).
Here, we answer this question.

\medskip\par\noindent{\bf Theorem~1. }{\em
There exists a non-constant analytic function $f\in \Sm$ satisfying
(1.4) and such that 
$\hbox{Re} f= 0$ a.e. on $\T$.
}

\medskip
On the other hand, condition (ii) can be really weakened:

\medskip\par\noindent{\bf Theorem~2. }{\em Let $f$ be an analytic function
in $\D$ satisfying conditions (i) and (iv), and let
$$
\liminf_{R\to\infty} R\int_R^\infty \frac{m_f(t)}{t}\, dt < \infty\,.
\eqno (1.5)
$$
Then $f$ is represented by the Schwarz integral (1.1) of a real measure
$\mu$, and}
$$
||\mu|| \le ||\hbox{Re} f||_1
+ C \liminf_{R\to\infty} R\int_R^\infty \frac{m_f(t)}{t}\, dt\,,
$$
{\em where $C$ is  a positive numerical constant.}

\medskip It is worth to note that after integration by parts condition
(1.5) can be written as
$$
\liminf_{R\to\infty} R\mathcal N\left(\frac{f}{R}\right) < +\infty\,,
$$
where
$$
\mathcal N(f) = \int_{\mathbb T} \log^+|f|\, dm
$$
is a ``norm'' in the Smirnov class $\mathcal N_+
$.

\bigskip\par\noindent{\bf \S 2. Proof of Theorem~1.}

\medskip\par\noindent 
We construct $f$ as a universal covering of $\C\setminus E$, where $E$ is 
a closed subset of $i\R_+$, by the unit
disc $\D$. We shall use some classical facts about universal coverings,
harmonic measures and Green functions \cite{N}.

Let $E=\bigcup_{n\ge 1} [ir_n, 2ir_n]$, where $r_1=1$, 
$r_{n+1}/r_n\to\infty$, and  let
$$
\o (t) = \o(0, E\cap \{|z|\ge t\}, \C\setminus E)
$$
be the harmonic measure of $E\cap \{z:\, |z|\ge t\}$ with respect to 
$\C\setminus E$ evaluated at the origin. We can choose the sequence 
$\{r_n\}$ increasing so fast that 
$$
\liminf_{t\to\infty} t \o (t) = 0\,.
\eqno (2.1)
$$
Indeed, let $h_n(z)$ be a bounded harmonic function in $\C\setminus 
\left( [i,2i]\bigcup [ir_n, i\infty)\,\right)$, vanishing on $[i,2i]$, and
having its boundary values equal identically to $1$ on $[ir_n,
i\infty)$. Evidently, $\lim_{n\to\infty} h_n(0) = 0$. Hence, on the
$n+1$-st step we can choose a value $r_{n+1}$ sufficiently large that 
$h_{n+1}(0) \le r_n^{-2}$. Then by the maximum principle
\begin{eqnarray*}
\liminf_{t\to\infty} t\o (t) &\le& \liminf_{n\to\infty} 2r_n\o (2r_n) \\
&=& 2 \liminf_{n\to\infty} r_n\o (r_{n+1}) \le 2 \liminf_{n\to\infty} r_n
h_{n+1}(0) = 0\,,
\end{eqnarray*}
proving (2.1). 

Now, let $f: \D \to \C\setminus E$ be the universal covering map
normalized by $f(0)=0$. 
Since $E$ has a positive capacity, $f$ is of bounded type in
$\D$, and therefore, has non-tangential boundary values a.e. on $\T$. By
the invariance of the harmonic measure,
$$
\int_\T (\varphi \circ f) (\z) \,\hbox{Re}\left(\frac{\z+z}{\z-z}\right) 
dm(\z) = 
\int_E \varphi (\eta)\, \o(f(z), d\eta, \C\setminus E)\,,
\eqno (2.2)
$$
where $\varphi$ is an arbitrary continuous function on $E$ (as the
boundary of $\C\setminus E$). In particular,
$$
\int_T (\varphi\circ f)\, dm = \int_E \varphi (\eta)
\, \o(0, d\eta, \C\setminus E)\,.
\eqno (2.2a)
$$
After a monotonic limit transitions, relations (2.2) and (2.2a) also hold
for semi-continuous functions, and therefore
$$
m_f(t) = \o(t)\,,
\qquad 0<t<\infty\,.
\eqno (2.3)
$$ 

Now, the function $f$ has pure imaginary boundary values a.e. on $\T$, and
$$
\liminf_{t\to\infty} tm_f(t) \stackrel{(2.3)}= \liminf_{t\to\infty} t\o(t)
\stackrel{(2.1)}= 0\,.
$$
It remains to observe that $f$ is of Smirnov's class $\Sm$. Indeed, 
representing a subharmonic function $\log|w|$ in $\C\setminus E$ as a sum
of the Green function and the Poisson integral, we have
$$
\log|w| = G_{\C\setminus E} (w,0) +
\int_E \log|\eta| \, \o(w, d\eta, \C\setminus E) + \alpha K_E(w)\,,
\eqno (2.4)
$$
where $G_{\C\setminus E} (w,0)$ is Green's function for $\C\setminus E$
with pole at $w=0$, $K_E(w)$ is a positive harmonic function in
$\C\setminus E$ vanishing on $E$ (so-called Martin function), and 
$\alpha\in \mathbb R$. 

First, we shall show that $\alpha \le 0$ (in fact, a minor
modification of the following argument shows that $\alpha =0$, though  
it is not needed for our purposes). Indeed, if $\alpha > 0$, then 
as follows from (2.4), the
function $K_E$ has a logarithmic growth at infinity. Setting
$K_E(z)=0$ for $z\in E$ we obtain a subharmonic function in $\C$ of
logarithmic growth, so that the set $E$ must be {\em thin} 
at infinity \cite{AH} and according to the Wiener criterion 
$$
\sum_{n=1}^\infty \frac{n}{\log \frac{2}{c(E_n)}} < \infty\,,
$$
where $E_n=E\cap \{z:\, 2^n\le |z| <2^{n+1}\}$ and  $c(\,.\,)$ is 
the logarithmic capacity (cf. \cite[Chapter~7]{Hayman} or
\cite[Chapter~5]{Ransford}). Since $E$ contains intervals $[ir_n, 2ir_n]$
and since the logarithmic capacity of the interval equals one quarter
of its length, we see that the series must diverge and therefore $\alpha
\le 0$. 

Omitting the positive terms on the RHS of (2.4), we get
$$
\log|w| \le
\int_E \log|\eta| \, \o(w, d\eta, \C\setminus E) \,.
$$
Setting here $w=f(z)$, $\eta = f(\z)$, and making use of (2.2), we obtain
$$
\log|f(z)| \le
\int_\T \log|f(\z)|\, 
\hbox{Re}\left(\frac{\z+z}{\z-z}\right) dm(\z)\,,
\qquad z\in\D\,, 
$$
which completes the proof. $\Box$

\bigskip\par\noindent{\bf \S 3. Proof of Theorem 2.}

\medskip\par\noindent We may assume
that $\hbox{Re} f = 0$ a.e. on $\T$. Otherwise, we decompose
$$
f(z) = f_1(z) + \int_\T \frac{\z+z}{\z-z} 
\left( \hbox{Re} f\right)(\z)\, dm(\z) = f_1 + f_2\,,
$$
where $\hbox{Re} f_1 = 0$ a.e. on $\T$, and
$\lim_{t\to\infty} tm_{f_2}(t) = 0$.

Following \cite{MOS}, \cite{MS1}, and \cite{MS2}, we
introduce a ``logarithmic determinant''
$$
u_f(w) := \int_\T \log|1-wf(\z)| \, dm(\z)\,.
\eqno (3.1)
$$
The function $u_f$ is subharmonic in $\C$, harmonic in the right and left
half-planes $\Pi_\pm$, and $u(0)=0$. Furthermore, 
$$
u_f(w) \le O(\log |w|)\,, 
\qquad w\to\infty\,,
$$
so that, $u_f$ is represented by the Poisson integrals in $\Pi_\pm$.
The function $u_f$ has a lower bound 
$$
u_f(it) \ge \log|1-it f(0)| \ge 0\,,
\qquad t\in \mathbb R.
$$
By the maximum principle applied to $u_f$ in $\Pi_\pm$, the function $u_f$
is non-negative everywhere in $\C$.

Our goal is to estimate from above the integral
$$
\I := \frac{1}{\pi} \int_{\mathbb R} \frac{u_f(it)}{t^2}\, dt\,.
$$
But first, we shall show that this integral controls 
the norm $||f||_{1,\infty}$ (cf. the proof of Theorem~2 in
\cite{MS1}). Indeed, estimating the Poisson
integrals in $\Pi_\pm$ we have
$$
u(re^{i\theta}) = \frac{r|\cos\theta|}{\pi} \int_\R 
\frac{u_f(it)}{t^2} \, \frac{t^2}{|it-re^{i\theta}|^2}\, dt
\le \frac{r}{|\cos\theta|} \, \I\,,
\qquad 0<r<\infty\,.
$$
In particular, $u_f(re^{i\theta})\le r\sqrt{2}\, \I$ within the angles 
$\{|\theta|\le \pi/4\}$ and $\{|\theta|\ge 3\pi/4 \}$.
Applying the maximum principle to the harmonic function
$u_f(z) - 2\I|\hbox{Im} z|$ within the complementary angles 
$\{\pi/4 \le |\theta| \le 3\pi/4\}$, we get
$$
M(r, u_f) := \max_{\theta\in [-\pi, \pi]} u_f(re^{i\theta}) 
\le 2\I r\,, 
\qquad 0<r<\infty\,.
$$
Let $\mu_f$ be the Riesz measure of the function $u_f$, and let
$\mu_f(r) = \mu_f\{|w|\le r\}$ be its counting function. Then
by the Jensen formula
$$
\mu_f(r) \le M(er, u_f) \le 2e\I r\,, \qquad  0<r<\infty\,.
$$ 
It remains to observe that 
$$
m_f(\tau) \stackrel{(3.1)}= \mu_f(1/\tau)\,,
\eqno (3.2)
$$
so that 
$$
||f||_{1,\infty} \le 2e \I\,. 
\eqno (3.3)
$$

Now, we estimate the integral $\I$ using an integral formula which has
been used previously in a similar situation (cf. the proof of Theorem~3 in
\cite{MOS}). For $|\theta| < \pi/2$ we have
$$
u_f(re^{i\theta}) = \frac{r\cos\theta}{\pi}
\int_\R \frac{u_f(it) dt}{|re^{i\theta} - it|^2}\,.
$$
Integrating this against $\cos\theta$, we get
\begin{eqnarray*}
\int_{-\pi/2}^{\pi/2} u_f(re^{i\theta}) \cos\theta\, d\theta
&=& \frac{r}{\pi} \int_\R u_f(it) \, dt 
\int_{-\pi/2}^{\pi/2} \frac{\cos^2\theta d\theta}{|re^{i\theta} - it|^2}
\\ \\ 
&=& \frac{r}{2} \int_\R u_f(it) \min\left( \frac{1}{t^2},
\frac{1}{r^2}\right) dt\,,
\end{eqnarray*}
since 
$$
\int_{-\pi/2}^{\pi/2} \frac{\cos^2\theta d\theta}{|re^{i\theta} - it|^2}
= \frac{\pi}{2} \min\left( \frac{1}{t^2}, \frac{1}{r^2}\right) \,.
$$
Using a similar relation in the left half-plane, we have
\begin{eqnarray*}
\int_\R u_f(it) \min\left( \frac{1}{t^2}, \frac{1}{r^2}\right) dt
&= & \frac{2}{r} \int_{-\pi}^\pi u_f(re^{i\theta}) |\cos\theta|\, d\theta
\\ \\ 
&< & \frac{2}{r} \int_{-\pi}^\pi u_f(re^{i\theta}) \, d\theta \\ \\
&= & \frac{4\pi}{r} \int_0^r \frac{\mu_f (s)}{s}\, ds\,.
\end{eqnarray*}
At last, making the monotonic limit
transition as $r\to 0$, and setting $R=1/r$, we get 
$$
\I \le  4\pi \liminf_{r\to 0} \frac{1}{r} \int_0^r \frac{\mu_f(s)}{s}\, ds 
\stackrel{(3.2)}\le 4\pi \liminf_{R\to\infty} R\int_R^\infty
\frac{m_f(\tau)}{\tau}\, d\tau\,.
\eqno (3.4)
$$
Combining inequalities (3.3) and (3.4) with Theorem~A1, we complete the
proof. $\Box$

\bigskip\par\noindent{\bf \S 4. Concluding remarks. }

\medskip\par\noindent{\bf 4.1 } The technique of logarithmic determinants
used in  the proof of Theorem~2 also allows us to prove another theorem of
Aleksandrov \cite{A1}, \cite{A2}:

\medskip\par\noindent{\bf Theorem A2. }{\em Let $f$ be an
analytic 
function in the unit disc satisfying condition (iv), let $f\in
\bigcap_{p<1} H^p$, 
and let 
$$
\liminf_{p\uparrow 1} (1-p)||f||_{H^p} < \infty\,.
$$
Then $f$ is represented by the Schwarz integral (1.1) of a real measure
$\mu$, and}
$$
||\mu|| \le ||\hbox{Re} f||_1 + C\liminf_{p\uparrow 1}
(1-p)||f||_{H^p}\,. 
\eqno (4.1)
$$

\noindent (In fact, Aleksandrov proved this estimate  
with $C=\frac{\pi}{2}$ on the RHS.)

Our proof follows the same lines as that of Theorem~2, only in the last
step we use other integral formulas to estimate the integral $\I$:
$$
\int_\R \log|1-it \lambda|\, \frac{dt}{|t|^{1+p}}
= |\lambda|^p\, \frac{\pi}{p} \cot \frac{\pi p}{2}\,,
\qquad \lambda\in \C\,, \quad 0<p<1\,,
$$
and
$$
\int_\R \frac{u_f(it)}{|t|^{1+p}}\, dt 
= \frac{\pi}{p} \cot \frac{\pi p}{2} \int_\T |f(\z)|^p\, dm(\z)\,,
\qquad 0<p<1\,.
$$

On the other hand, Theorem~A2 follows directly from our Theorem~2 since 
the $p$-th power of the $L^p$ norm can be written as the integral
$$
p\int_0^\infty t^p\, \frac{m_f(t)}{t}\, dt
= p^2 \int_0^\infty t^{p-1} 
\left( \int_t^\infty \frac{m_f(s)}{s}\, ds \right) dt\,.
$$
This yields that if 
$$
\lim_{t\to \infty} t \int_t^\infty \frac{m_f(s)}{s}\, ds = \infty\,,
$$
then $\lim_{p\uparrow 1} (1-p)||f||_p = \infty$.

We also note that it is not difficult to construct a 
function $h$ on
$[-\pi,\pi]$ such that 
$$
\liminf_{t\to\infty} t\int_t^\infty
\frac{m_h(s)}{s}\, ds = 0
$$
while $h$ does not belong to any of the $L^p$
spaces for $p>0$. This shows that the assumptions of Theorem~2 are really
weaker than those of Theorem~A2.

\smallskip\par\noindent{\bf 4.2 } The counter-example provided by
Theorem~1 is related to the logarithmic determinants by the
following formula (in the notation of $\S 2$):
$$
\int_\T \log|1-wf(\z)|\, dm(\z) \stackrel{(2.2a)}= 
G_{\C\setminus E^*}(w,\infty)\,,
$$
where $E^*=\{\lambda:\, 1/\lambda\in E\}$.

\smallskip\par\noindent{\bf 4.3 }
Theorems~1 and 2 can be easily reformulated for analytic
functions represented by Cauchy-type integrals
$$
f(z) = \int_\T \frac{d\mu(\zeta)}{\zeta-z}
$$
of complex-valued measures $\mu$ of finite variations. In this case,
condition (iv) must be replaced by 

\smallskip\par\noindent{\bf (iv')} $Tf\in L^1(\T)$, where 
$Tf(\zeta) = \lim_{r\uparrow 1} \left(f(r\zeta) - f(r^{-1}\zeta)
\right)$.

\smallskip\par\noindent
We leave the details to the reader (cf. \cite{A1}, \cite{A2}).

\bigskip\par\noindent{\bf Acknowledgements.} I thank A.~Aleksandrov,
I.~Ostrovskii and P.~Yuditskii for their numerous useful remarks and
suggestions.

\bigskip\par\noindent{\em School of Mathematical
Sciences, \newline Tel-Aviv University, \newline
\noindent Ramat-Aviv, 69978, \newline Israel

\smallskip\par\noindent sodin@post.tau.ac.il}

\end{document}